# STABILIZATION OF UNDERACTUATED MECHANICAL SYSTEMS WITH TIME-VARYING UNCERTAINTY USING ADAPTIVE FUZZY SLIDING MODE

Mohammad Mahdi Azimi, Hamid Reza Koofigar*, Mehdi Edrisi

Faculty of Engineering, University of Isfahan, Isfahan 81746-73441, Iran

## ABSTRACT

**In this paper, the control problem of underactuated systems in the presence of external disturbances and time varying model uncertainties is considered. Ensuring the stability and robust performance against uncertainties and disturbances are the main problems in such systems. An adaptive fuzzy sliding mode controller (AFSMC) is proposed to solve the problem, satisfying the robustness properties against the perturbations. The designed controller can be applied to a wide class of underactuated systems with fewer restrictions, compared with many pervious works. To evaluate the performance of proposed algorithm, it is applied to an underactuated inverted pendulum and a Translational Oscillator/Rotational Actuator (TORA) system with two degrees of freedom. The simulation results are compared with a conventional method to demonstrate the effectiveness of the proposed strategy. The results show that the proposed controller guarantees the stability and robustness against time varying uncertainties and disturbances.**

## 1. INTRODUCTION

Control of underactuated systems is especially important, due to extensive research areas and many applications in robotics, aerospace, naval vessels and submarines. In these systems, the number of actuators is less than the outputs; therefore, some degrees of freedom are not excited directly by the actuators. Reasons such as the dynamics of systems, systems design, fault in actuators and some conditions for model reduction can make the system underactuated. In the past decades, considerable attentions have been paid to control problem of underactuated systems, as studied by [1, 2]. More recently, [3] have developed sliding mode controller, as an effective strategy against uncertainties, in such important applications as self-balancing robots, mobile robots and submarines [4]. An approach to design a sliding mode controller for a specific class of underactuated mechanical systems in cascade form has been presented by [5], and based on a novel transformation developed by [6, 7]. In fact, the method provides decoupling, using a systematic approach to transform a class of underactuated mechanical systems into a cascade form. Stabilizing underactuated systems in the presence of matched and mismatched uncertainties based on lyapunov redesign method, is introduced by

[8]. Also, [9] have presented an adaptive control method, based on the transformation introduced by [7].

Sliding mode controller has been presented for a class of second-order underactuated systems which does not guarantee the asymptotic stability, [10]. Removing such drawback, an adaptive backstepping algorithm, using wavelet network, has been proposed with some design complexities, [11].

In recent researches dynamic system with time-varying uncertainties considered in many papers. In [12] controller designed for linear system with time varying uncertainties. Also, sliding control for nonlinear system containing time-varying uncertainties with unknown bounds is studied [13]. In some of papers, the proposed controller is implemented to a single-link flexible robot manipulator [14] and flexible ball screw drives [15].

In this paper, removing some conservative assumptions of the pervious works, an adaptive fuzzy control algorithm is proposed to ensure robust asymptotic stability for a general class of underactuated systems, with possessing simplicity and universality properties, as two main benefits of the method. In addition to mathematical analysis of the proposed scheme, a comparison study is also made on two benchmark problems in mechanical engineering, by demonstrating some simulation results.

This paper is organized as follows: In section 2, the stages of problem formulation of underactuated system and implied time varying uncertainties and disturbance is described. In section 3, controller design procedure is described. In section 4, two numerical examples are given and the results and discussion is presented. Finally, section 5 concludes the paper.

## 2. PROBLEM FORMULATION

As a preliminary step to controller design procedure, a class of nonlinear underactuated systems, perturbed by uncertainties and external disturbances, introduced. Consider the underactuated systems described by

$$M(\theta)\ddot{\theta} + C(\theta,\dot{\theta})\dot{\theta} + G(\theta) = T \qquad (1)$$

where

$$\theta = \begin{pmatrix} \theta_1 \\ \theta_2 \end{pmatrix}, T = \begin{pmatrix} 1 \\ 0 \end{pmatrix}\tau$$

and

$$M(\theta) = \begin{bmatrix} m_{11}(\theta) & m_{12}(\theta) \\ m_{21}(\theta) & m_{22}(\theta) \end{bmatrix}$$

$$C(\theta,\dot{\theta}) = \begin{bmatrix} c_{11}(\theta,\dot{\theta}) & c_{12}(\theta,\dot{\theta}) \\ c_{21}(\theta,\dot{\theta}) & c_{22}(\theta,\dot{\theta}) \end{bmatrix}$$

$$G(\theta) = \begin{bmatrix} g_1(\theta_1) \\ g_2(\theta_2) \end{bmatrix}$$

are the inertia matrix, the matrix of the coriolis and centrifugal forces, and the vector of gravity forces respectively.

In order to synthesis the adaptive backstepping controller, the affine form of dynamical equation (1) is obtained. To this end, take $x_1 = \theta_1$, $x_2 = \dot{\theta}_1$, $x_3 = \theta_2$, $x_4 = \dot{\theta}_2$, and rewrite the equations in presence of matched and mismatched uncertainties and external disturbance as

$$\begin{aligned} \dot{x}_1 &= x_2 \\ \dot{x}_2 &= f_1(\mathbf{x}) + g_1(\mathbf{x})u + d_1(t) \\ \dot{x}_3 &= x_4 \\ \dot{x}_4 &= f_2(\mathbf{x}) + g_2(\mathbf{x})u + d_2(t) \end{aligned} \qquad (2)$$

where $\mathbf{x} = \begin{bmatrix} x_1, x_2, x_3, x_4 \end{bmatrix}^T$ is vector of states, $f_1(\mathbf{x})$, $f_2(\mathbf{x})$, $g_1(\mathbf{x})$ and $g_2(\mathbf{x})$ are unknown nonlinear functions, $u$ is control input and $d_1(t)$, $d_2(t)$ are external disturbances, i.e., $|d_1(t)| \leq D_1(t)$, $|d_2(t)| \leq D_2(t)$. Defining $y = [x_1, x_3]^T$ as outputs vector, the system can be considered as an underactuated system, and input is $u$. As shown in equation (1), the number of inputs is less than outputs, therefore the system can be considered as single-input multi-output (SIMO) underactuated system. To design a conventional sliding mode controller, the system dynamical equations should be in a canonical form, [16]. Hence, designing the sliding controller for the system under study (2), is not straightforward. Moreover, actuated and non-actuated degrees of freedom are coupled together, so it is must be considered in controller design.

The control problem is to design an Adaptive Fuzzy Sliding Mode Controller (AFSMC) for underactuated system (2) such that the boundedness of all closed loop signals and asymptotic stability are guaranteed.

## 3.     CONTROLLER DESIGN

In order to design the AFSMC, two sliding surface are first defined as

$$s_1 = c_1 x_1 + x_2 \qquad (3)$$

and

$$s_2 = c_2 x_3 + x_4 \qquad (4)$$

where $c_1$ and $c_2$ are two positive constants. According to the principle of sliding mode controller design, the control input can be selected as follows

$$u_1 = \hat{u}_1 - K_1 \mathrm{sat}(s_1 g_1(\mathbf{x}) / \Phi_1), \ K_1 > D_1 / |g_1(\mathbf{x})| \qquad (5)$$

with

$$\hat{u}_1 = \frac{-c_1 x_2 - f_1(\mathbf{x})}{g_1(\mathbf{x})} \qquad (6)$$

or

$$u_2 = \hat{u}_2 - K_2 \mathrm{sat}(s_2 g_2(\mathbf{x}) / \Phi_2) \\ , \ K_2 > D_2 / |g_2(\mathbf{x})| \qquad (7)$$

with

$$\hat{u}_2 = \frac{-c_2 x_4 - f_2(\mathbf{x})}{g_2(\mathbf{x})} \qquad (8)$$

Obviously, if $u = u_1$ used to control the system, then the states $x_1$ and $x_2$ are asymptotically converged to zero along the surface $s_1$, and if $u = u_2$ is used, the convergence is obtained for the states $x_3$, $x_4$ along the surface $s_2$. Thus, in order to control of all states, only by the one control input $u_1$, the sliding surface $s_1$ can be modified to

$$s_1 = c_1(x_1 - z) + x_2 \qquad (9)$$

where $z$ is the proportional to $s_2$, transferred from $s_2$. Equation (9) implies that, the control input $u_1$ is changed in such a way that $x_1$ and $x_2$ converges to $x_1 = z$ and $x_2 = 0$, [17].

Since, the control input used to control the system, should be guaranteed to be bounded, suppose

$$|z| < z_U, \ 0 < z_U < 1 \qquad (10)$$

which $z$ has been chosen as follows

$$z = \mathrm{sat}(s_2 / \Phi_z) \cdot z_U \qquad (11)$$

where, the $\Phi_z$ is boundary layer for $s_2$, used to smoothen $z$. The value of $\Phi_z$, puts the sliding surface $s_2$ in a suitable range of $x_1$, and sat($\cdot$) has been defined as follows

$$\mathrm{sat}(\varphi) = \begin{cases} \mathrm{sgn}(\varphi), \text{if } |\varphi| \geq 1 \\ \varphi, \text{if } |\varphi| < 1 \end{cases} \qquad (12)$$

To generate the control input for control purpose, the equivalent term in (5) is considered in a general form as

$$\begin{aligned} \hat{u}_1 = u = \\ \frac{1}{\hat{g}_1(x \mid \theta_g)} [-c_1 x_2 - \hat{f}_1(x \mid \theta_f) - \\ K_p \mathrm{sat}(s_1 \hat{g}_1(x \mid \theta_g) / \Phi_1)] \end{aligned} \qquad (13)$$

with $\hat{g}_1(\mathbf{x} \mid \theta_g) > 0$. The estimation of unknown functions $f_1(\mathbf{x})$ and $g_1(\mathbf{x})$ are respectively denoted by

$$\hat{f}_1(\mathbf{x} \mid \theta_f) = \theta_f^T \xi(\mathbf{x}) \qquad (14)$$

$$\hat{g}_1(\mathbf{x} \mid \theta_g) = \theta_g^T \eta(\mathbf{x}) \qquad (15)$$

calculated by an adaptive fuzzy technique.

***Theorem 1***. Consider the nonlinear underactuated system (2). If the control input (13) is used, and functions $\hat{f}_1(x \mid \theta_f)$ and $\hat{g}_1(x \mid \theta_g)$ are estimated by (14) and (15) with parameter adjustment mechanisms

$$\dot{\theta}_f = \gamma_1 s_1 \xi(\mathbf{x}) \qquad (16)$$

$$\dot{\theta}_g = \gamma_2 s_1 \eta(\mathbf{x}) u \qquad (17)$$

then all the closed-loop signals remain bounded and the tracking error is converge to zero asymptotically.

***Proof:*** The optimized parameter vectors can be defined as

$$\theta_f^* = \arg \min_{\theta_f \in \Omega_f} [sup \mid \hat{f}(\mathbf{x} \mid \theta_f) - f(\mathbf{x}) \mid] \qquad (18)$$

$$\theta_g^* = \arg \min_{\theta_g \in \Omega_g} [sup \mid \hat{g}(\mathbf{x} \mid \theta_g) - g(\mathbf{x}) \mid] \qquad (19)$$

where $\Omega_f$ and $\Omega_g$ are constraint sets for $\theta_f$ and $\theta_g$ as

$$\Omega_f = \{\theta_f \in R^n \mid \left|\theta_f\right| \le M_f\} \qquad (20)$$

$$\Omega_g = \{\theta_g \in R^n \mid 0 < \varepsilon \le \left|\theta_g\right| \le M_g\} \qquad (21)$$

where $M_f$, $\varepsilon$, and $M_g$ are some predefined bounds. Now, the minimum approximation error can be defined as

$$w = [\hat{f}(\mathbf{x} \mid \theta_f^*) - f(\mathbf{x}, t)] + [\hat{g}(\mathbf{x} \mid \theta_g^*) - g(\mathbf{x}, t)] u \qquad (22)$$

It is assumed that the fuzzy parameters $\theta_f$ and $\theta_g$ never reach to the defined bounds, [18].

Now, a lyapunov function candidate is choose as

$$V = \frac{1}{2} s_1^2 + \frac{1}{2\gamma_1} \tilde{\theta}_f^T \tilde{\theta}_f + \frac{1}{2\gamma_2} \tilde{\theta}_g^T \tilde{\theta}_g \qquad (23)$$

where $\tilde{\theta}_f = \theta_f^* - \theta_f$ and $\tilde{\theta}_g = \theta_g^* - \theta_g$ denote the parameter estimation errors. The time derivate of lyapunov function (23) is calculated as

$$\dot{V} = s_1 \dot{s}_1 + \frac{1}{\gamma_1} \tilde{\theta}_f^T \dot{\tilde{\theta}}_f + \frac{1}{\gamma_2} \tilde{\theta}_g^T \dot{\tilde{\theta}}_g \qquad (24)$$

On the other hand, by differentiating (9) and using definition (22), one can obtain

$$\begin{aligned}
\dot{s}_1 &= c_1 x_2 - c_1 \dot{z} + \hat{f}_1(\mathbf{x} \mid \theta_f^*) + \hat{g}_1(\mathbf{x} \mid \theta_g^*) u + w + d_1(t) \\
&= c_1 x_2 - c_1 z + \hat{f}_1(\mathbf{x} \mid \theta_f^*) - \hat{f}_1(\mathbf{x} \mid \theta_f) \\
&\quad + \left(\hat{g}_1(\mathbf{x} \mid \theta_g^*) - \hat{g}_1(\mathbf{x} \mid \theta_g)\right) u \\
&\quad + \hat{f}_1(\mathbf{x} \mid \theta_f) + \hat{g}_1(\mathbf{x} \mid \theta_g) u + w + d_1(t)
\end{aligned} \qquad (25)$$

Incorporating adaptation laws (16) and (17) and control input (13) into (25) yield

$$\dot{s}_1 = -c_1 \dot{z} + \tilde{\theta}_f^T \xi(\mathbf{x}) + \tilde{\theta}_g^T \eta(\mathbf{x}) u - p(\mathbf{x}) + w + d_1(t) \qquad (26)$$

where $p(\mathbf{x}) = K_p \text{sat}(s_1 \hat{g}_1(\mathbf{x} \mid \theta_g) / \Phi_1)$. Substituting (26) into (24) yields

$$\begin{aligned}
\dot{V} &= s_1 \tilde{\theta}_f^T \xi(\mathbf{x}) + \frac{1}{\gamma_1} \tilde{\theta}_f^T \dot{\tilde{\theta}}_f \\
&\quad + s_1 \tilde{\theta}_g^T \eta(\mathbf{x}) u + \frac{1}{\gamma_2} \tilde{\theta}_g^T \dot{\tilde{\theta}}_g \\
&\quad - s_1 c_1 \dot{z} - s_1 p(\mathbf{x}) + s_1 (w + d_1(t))
\end{aligned} \qquad (27)$$

By $\dot{\tilde{\theta}}_f = -\dot{\theta}_f$ and $\dot{\tilde{\theta}}_g = -\dot{\theta}_g$, and taking the update laws (16) and (17) in to account, one obtains

$$\dot{V} \le -s_1 c_1 \dot{z} - s_1 K_p \text{sat}(s_1 \hat{g}_1(\mathbf{x} \mid \theta_g) / \Phi_1) + s_1 w + D_1 |s_1| \qquad (28)$$

By choosing $K_p = D_1 + \eta_1$, the inequality (28) takes the form

$$\dot{V} \le -s_1 c_1 \dot{z} - K_p s_1 \text{sgn}(s_1) + s_1 w + D_1 |s_1| \qquad (29)$$

Now, set $K_p = D_1 + \eta_1$, $\eta_1 > 0$ to get

$$\dot{V} \le -s_1 h - \eta_1 |s_1| + s_1 w \qquad (30)$$

where $h = c_1 \dot{z}$. By definitions (9) and (11) and sufficiently large choice of $\Phi_z$, $\dot{z} \to 0$ as $s_1 \to 0$. Therefore, all signals in the system are bounded. Obviously, if $x_1(0)$ is bounded, then $x_1(t)$ for all $t$ is bounded, therefore all states of system $\mathbf{x}$ will be bounded. To complete the proof, we need proving that $\dot{s}_1 \to 0$ as $t \to \infty$.

Assume that $|s_1| \le \eta_s$, then (26) can be rewritten as

$$\dot{V} \le |s_1| w - |s_1| \eta_1 \le \eta_s |w| - |s_1| \eta_1 \qquad (31)$$

Integrating both sides of (31) yields

$$\int_0^\tau |s_1| d\tau \le \frac{1}{\eta} (|V(0)| + |V(t)|) + \frac{\eta_s}{\eta_1} \int_0^\tau |w| d\tau \qquad (32)$$

Then, we have $s_1 \in L_1$ and every term in (23) is bounded. Hence, $s_1, \dot{s}_1 \in L_\infty$ use for Barbalat's lemma, [16]. We have $s_1(t) \to 0$ as $t \to \infty$ therefore $z \to 0$ and $s_2 \to 0$, the system is stable and the error will asymptotically converge to zero.

∎

**Remark I.** For design a conventional sliding mode controller, the system dynamic equations must be in canonical form [16]. This assumption has been removed here for underactuated system (1).

**Remark II.** Although there exist four unknown functions in the system under consideration, only two function approximations are needed for control effort calculation. This facilitates implementation of the proposed control algorithm.

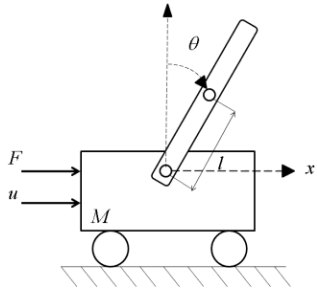

**Fig. 1 Inverted Pendulum**

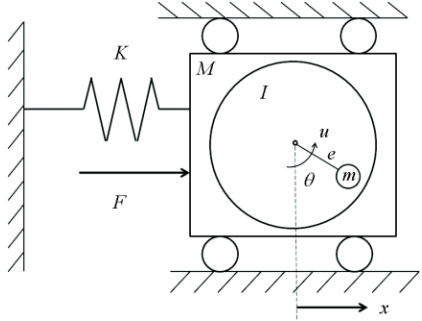

**Fig. 2 Translational Oscillator/Rotational Actuator (TORA)**

# 4. RESULTS AND DISCUSSION

To evaluate the performance of the proposed strategy, it is applied to an inverted pendulum and a Translational Oscillator/Rotational Actuator (TORA) system, as two underactuated benchmark systems which shown in Figs. 1 and 2.

## 4.1. Application to Inverted Pendulum

The dynamic equation of an inverted pendulum can be described in the form of (2) with

$$f_1 = \frac{m_t g \sin x_1 - m_p L \sin x_1 \cos x_1 x_2^2}{L(\frac{4}{3}m_t - m_p \cos^2 x_1)} \qquad (29)$$

$$g_1 = \frac{\cos x_1}{L(\frac{4}{3}m_t - m_p \cos^2 x_1)} \qquad (30)$$

$$f_2 = \frac{-\frac{4}{3}m_p L x_2^2 \sin x_1 + m_p g \sin x_1 \cos x_1}{\frac{4}{3}m_t - m_p \cos^2 x_1} \qquad (31)$$

$$g_2 = \frac{4}{3(\frac{4}{3}m_t - m_p \cos^2 x_1)} \qquad (32)$$

where $x_1(rad)$ denotes the angle of the pole with respect to vertical axis, $x_3(m)$ is the cart displacement, $m_c$ represents the mass of the cart, $m_p$ is the pole mass and $m_t = m_c + m_p$. Moreover, $f_1(\mathbf{x})$ and $g_1(\mathbf{x})$ is estimated by using adaptive fuzzy. Let the initial states vector be $\mathbf{x} = [\pi/15, 0, 0, 0]^T$ and desired states be $\mathbf{x}_d = [0, 0, 1, 0]^T$. In other words, the control objective is the cart displacement from $0(m)$ to $1(m)$ and pole angle from $\pi/15(rad)$ to $0(rad)$. The controller parameter are selected as $c_1 = 5$, $c_2 = 0.5$, $\phi_1 = 5$, $\phi_z = 15$, $K_p = 10$, $z_v = 0.945$ and positive constants in adaptation law as $\gamma_1 = 100$, $\gamma_2 = 2$. Three cases are considered here for simulation. In the first case, the nominal dynamical model of inverted pendulum is adopted. The second case takes into account a significant variation in system masses. In the third case, the adaptive fuzzy algorithm is applied in the presence of external disturbances.

*Case I.* Control of nominal inverted pendulum. Assume the pole angle is restricted to $z_3 \in [-\pi/6, \pi/6]$. The results of this section are shown in Fig 3.

In this case, there is no uncertainty in model of system and input, and controller performance is evaluated in an ideal situation. The pole angle becomes zero after a short time and cart placed in desired location. The control input is bounded and small enough, and achieved without causing saturation in actuators.

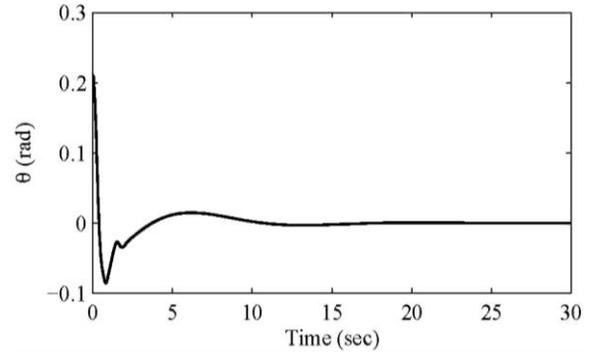

(a)

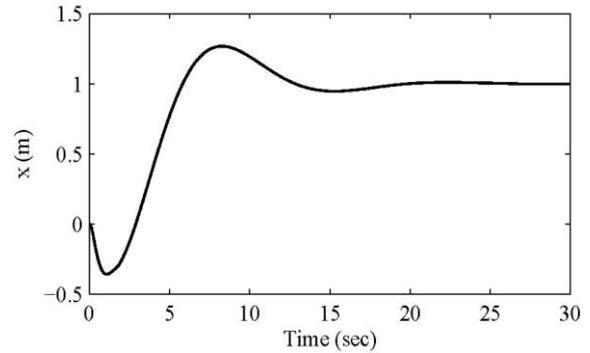

(b)

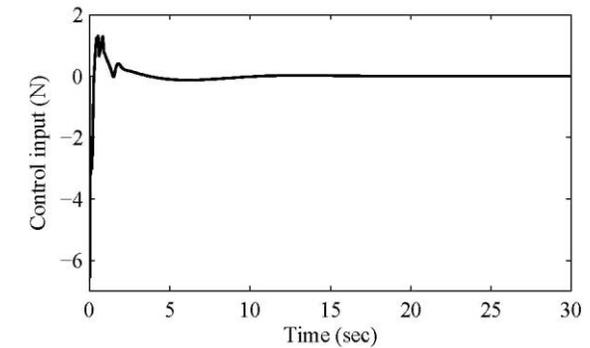

(c)

**Fig. 3 The simulation results for Inverted Pendulum in Case I, a) Pole angle, b) Cart displacement, c) Control input**

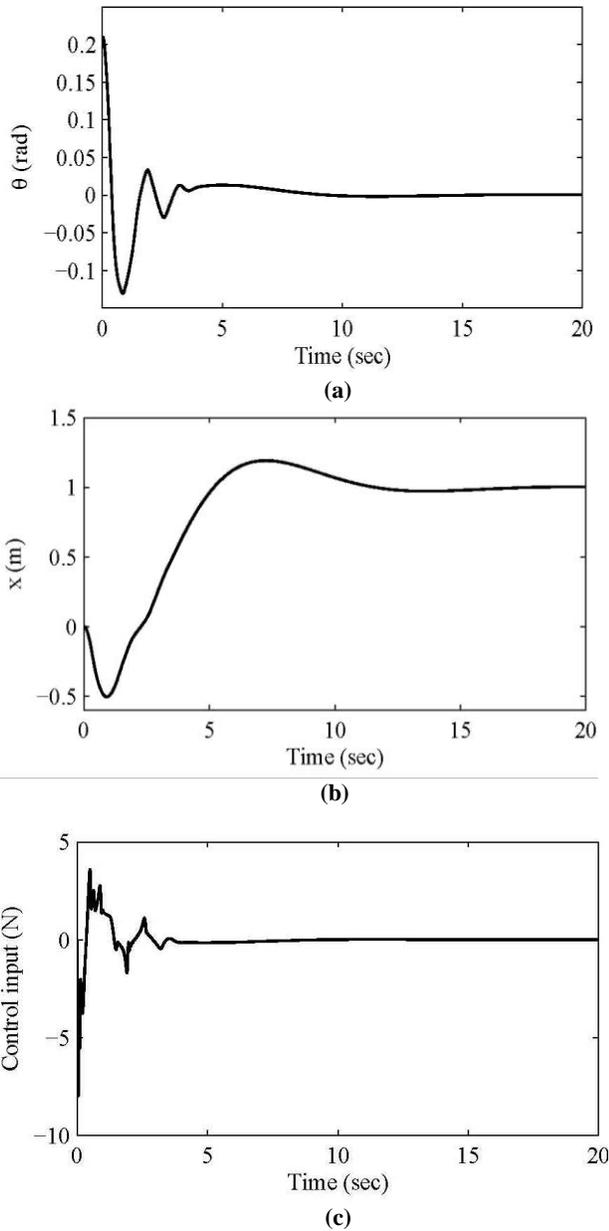

**Fig. 4 The simulation results for Inverted Pendulum in Case II, a) Pole angle, b) Cart displacement, c) Control input.**

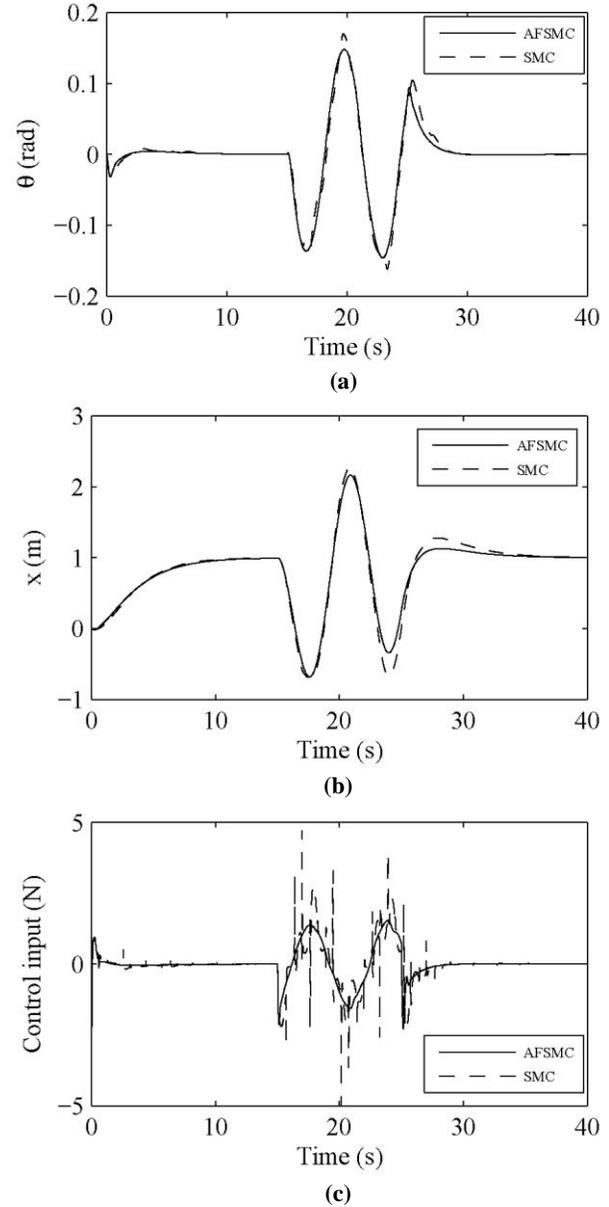

**Fig. 5 The simulation results for Inverted Pendulum in Case III, a) Pole angle, b) Cart displacement, c) Control input.**

*Case II.* The masses of the cart and pole are allowed to vary as $m_p = 0.1 + 0.05sin(t)$ (kg) and $m_c = 1 + 0.5sin(t)$ (kg). The results of this section are shown in Fig. 4.

The simulation results show that the considerable change in pole and cart masses has no significant effect on the closed loop performance and the desired performance of the controller is confirmed. Also the control input is acceptable.

*Case III.* The sinusoidal disturbance $F = 2.5\cos 5t$ is applied to the cart during $t \in [15, 25]$ to evaluate the robustness property of the strategy.

In this case, the external force in a specific interval of time is also taken, in addition to model uncertainty in case II. The simulation results in Fig. 5 show that the changes in pole angle and cart position are small, and the controller is still able to bring the states regulate at desired states.

Also, AFSMC simulation results are compared with decoupled sliding mode controller (SMC) [17] in Fig. 5. By applying the external disturbance to inverted pendulum, the AFSMC has more ability to attenuate the disturbance effect at the output with respect to SMC. Also, sliding control input has bigger swings and variations, compared with AFSMC. Hence, one can conclude the AFSMC has more robustness against SMC, and even with time-varying masses and various uncertainties the controller is still able to regulate the states to the desired trajectories.

## 4.1. Application to TORA

The dynamic equation for a TORA [9], is well known to be of form (2) with

$$f_1(\mathbf{x}) = \frac{-x_1 + \varepsilon x_4^2 \sin(x_3)}{1 - \varepsilon^2 \cos^2 x_3} \tag{33}$$

$$g_1(x_3) = -\frac{\varepsilon \cos x_3}{1 - \varepsilon^2 \cos^2 x_3} \tag{34}$$

$$f_2(\mathbf{x}) = \frac{\varepsilon \cos x_3 (x_1 - \varepsilon x_4^2 \sin x_3)}{1 - \varepsilon^2 \cos^2 x_3} \tag{35}$$

$$g_2(x_3) = \frac{1}{1 - \varepsilon^2 \cos^2 x_3} \tag{36}$$

where $x_1(m)$ is the linear displacement of platform, $x_3(rad)$ is the rotor displacement, The parameter $\varepsilon$ is defined as [8].

$$\varepsilon = \frac{me}{\sqrt{(I + me^2)(M + m)}} \tag{37}$$

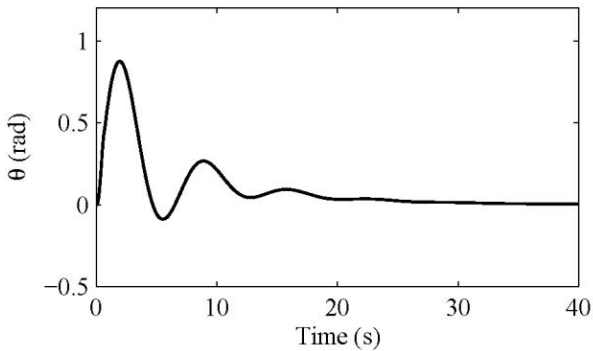

**(a)**

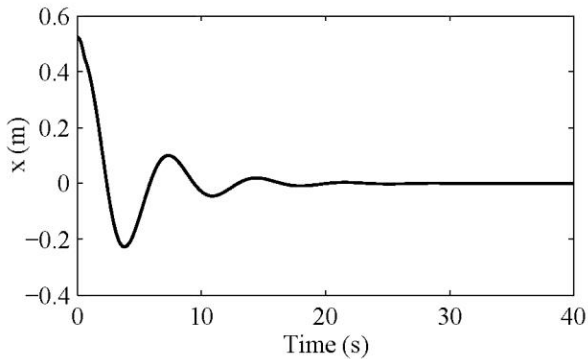

**(b)**

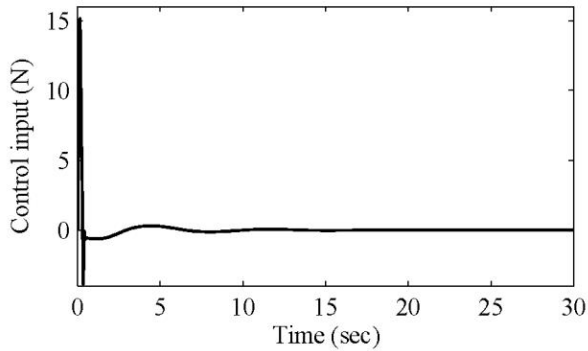

**(c)**

**Fig. 6 The simulation results for TORA in Case I, a) Angle, b) Position, c) Control input.**

where $m(kg)$ and $M(kg)$ are, respectively, the mass of the rotor and platform. Moreover $e(kg)$ is the eccentricity and $I(kgm^2)$ is the rotary inertia. Assume that $f_1(\mathbf{x})$, $g_1(x_3)$ and $f_2(\mathbf{x})$, $g_2(x_3)$ contain uncertainties and a stabilizing controller is needed to be designed.

The proposed controller is employed with the initial condition $\mathbf{x}(0) = \begin{bmatrix} \pi/6, 0, 0, 0 \end{bmatrix}$, and the desired state $\mathbf{x}_d(0) = \begin{bmatrix} 0, 0, 0, 0 \end{bmatrix}$. To avoid possible singularity in the system dynamics, we choose $\varepsilon = 0.33$ which implies that the actual values for the plant parameters are $m = 0.5 kg$, $M = 2 kg$, $I = 0.1 kgm^2$ and $e = 0.5 m$ [9]. The controller parameter are selected as $c_1 = 5$, $c_2 = 0.1$, $\phi_1 = 10$, $\phi_2 = 1$, $K_p = 120$, $z_v = 0.6$ and $\gamma_1 = 2$, $\gamma_2 = 100$. Two simulation cases are presented here. The first case is the control of an uncertain TORA. In the second case, there exists a significant variation in system masses in presence of external disturbance.

*Case I.* the masses of the cart and pole are allowed to vary with significant amount as $m = 0.5 + 0.1 sin(t)$ (kg) and $M = 2 + 0.4 sin(t)$ (kg). The results of this section are shown in Fig. 6.

In this case, first part of uncertainties is considered as model uncertainties. It is assumed that part of mass of pole and cart significantly are changed and controller performance is evaluated in this situation. As can be seen significant change in pole and cart masses has no significant effect on controller performance and proper performance of the controller is confirmed.

*Case II.* A sinusoidal disturbance $F = 0.25 \cos t$ is applied to the cart during $t \in \begin{bmatrix} 30, 40 \end{bmatrix}$ to test robustness of the strategy. The results of this section are shown in Fig. 7.

In this case, another uncertainty added to system. Input uncertainty added to system as external force in specific interval of time, also assumed that model uncertainties is there like last case, and controller performance is evaluated in this situation. It can be seen that the change in angle and position in desired state is small. Also change in control input can be seen in time that external force is applied and this is corresponds to input uncertainty.

The AFSMC simulation results are also compared with decoupled SMC [17] in Fig. 7, which shows faster convergence of system outputs to the zero.

## 5. CONCLUSION

In this paper, an adaptive fuzzy sliding mode controller has been proposed for a special class of underactuated systems. Compared some previous works, such control strategy possesses the simplicity and universality properties and ensures robust tracking performance, despite the perturbations. By introducing a practical underactuated system, the AFSMC has been implemented to two typical case studies including the Inverted pendulum and TORA. From a comparison point of view, the simulation results have been presented and discussed. The speed of convergence and lower control effort are two main benefits of the proposed scheme, compared with those of a sliding mode controller.

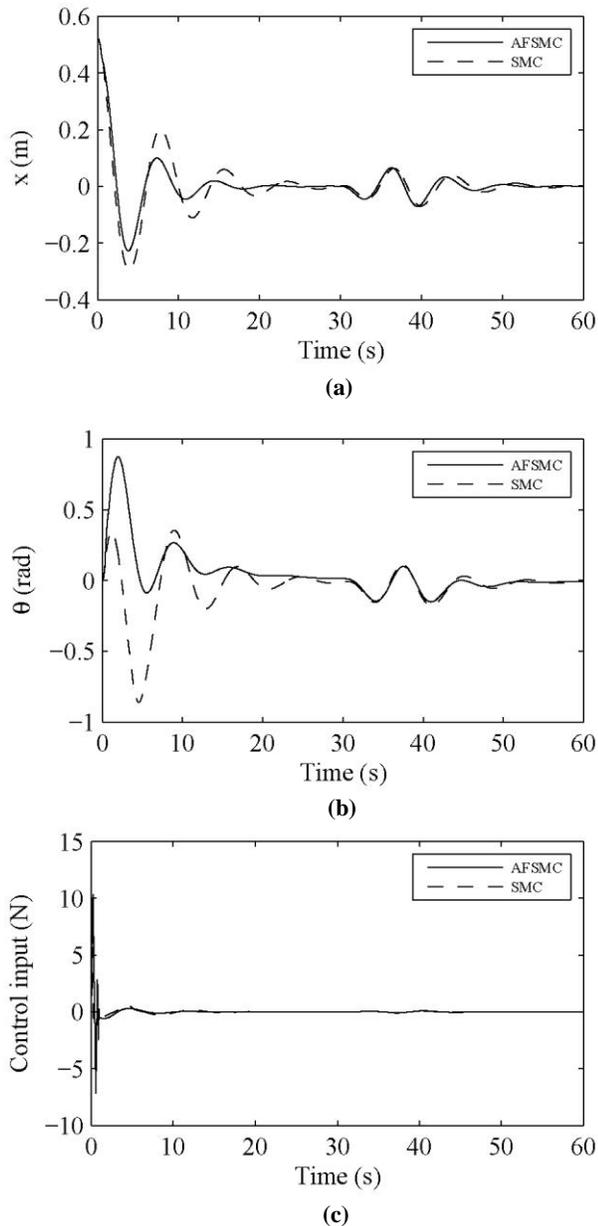

**Fig. 7 The simulation results for TORA in Case II, a) Angle, b) Position, c) Control input.**

*Mohammad Mahdi Azimi* is a PhD student in Control Engineering at Shahid Beheshti University, Tehran, Iran. He received the M.Sc. in Control Engineering from University of Isfahan, Isfahan, Iran, in 2013. His current research interests are nonlinear control, robust control, adaptive control, intelligent control and robotic systems.

*Hamid Reza Koofigar* received his BSc in Electronics, MSc in Control Engineering, and Ph.D. in Electrical Engineering from Isfahan University of Technology, Iran, in 2003, 2005, and 2009 respectively. Since Jan. 2010 he is with the Department of Electrical Engineering-Control Group, University of Isfahan, as an assistant professor. His current research interests include robust control, adaptive nonlinear control, intelligent control and robotics.

*Mehdi Edrisi* received his BSc in Electrical Engineering from Isfahan University of Technology, Iran, in 1990, MSc in Control Engineering from University of New South Wales, Australia, in 1994, and Ph.D. in Electrical Engineering from University of South Australia, Australia, in 2000. His current research interests include Control, Microprocessor Programming, Computer Simulation, EMC, Computer network design, Internet Technology.